\documentclass[letterpaper, 10 pt, journal, twoside]{ieeetran}  %

\usepackage[english]{babel}
\usepackage[T1]{fontenc}
\usepackage{graphicx} 
\graphicspath{{figures/}} 
\usepackage{amsmath}
\usepackage{amsfonts}
\usepackage{amssymb}
\usepackage{arydshln}
\usepackage{bm}
\usepackage{cite}
\usepackage[linesnumbered, noline, ruled]{algorithm2e}
\usepackage{etoolbox}
\makeatletter
\patchcmd{\@algocf@start}
{-1.5em}
{0pt}
{}{}
\makeatother
\usepackage{booktabs}

\usepackage{enumitem}
\usepackage{xcolor}

\usepackage{relsize}
\usepackage{exscale}
\usepackage[absolute]{textpos} 
\usepackage{hyperref}

\renewcommand{\r}[1]{(\ref{#1})}

\newcommand{\R}{\mathbb{R}}         
\newcommand{\F}{\mathbb{F}}
\renewcommand{\P}{\mathbb{P}}
\newcommand{\I}{\mathbb{I}}

\newcommand{\Fc}{\mathcal{F}}

\newcommand{\Ic}{\mathcal{I}}

\newcommand{\Ab}{\mathbf{A}}

\newcommand{\Cb}{\mathbf{C}}

\newcommand{\cl}{\prec}
\newcommand{\cg}{\succ}
\newcommand{\Delf}{\mathbf{\Delta}}

\newcommand{\diag}{\mathrm{diag}}

\newcommand{\Del}{\Delta}
\newcommand{\del}{\delta}
\newcommand{\ga}{\gamma}
\newcommand{\eps}{\varepsilon}

\DeclareMathOperator*{\argmin}{arg\,min}

\newcommand{\mat}[2]{\left(\begin{array}{@{}#1@{}}#2\end{array}\right)} 
\newcommand{\smat}[1]{\left(\begin{smallmatrix}#1\end{smallmatrix}\right)}

\newenvironment{red_test}{\color{red}}{} 

\newenvironment{blue_test}{\color{blue}}{} 

\renewcommand{\t}{\tilde}

\newtheorem{theorem}{Theorem}

\newtheorem{remark}[theorem]{Remark}

\def\BibTeX{{\rm B\kern-.05em{\sc i\kern-.025em b}\kern-.08em
		T\kern-.1667em\lower.7ex\hbox{E}\kern-.125emX}}
\title{Controller Design via Experimental Exploration with Robustness Guarantees
}

\author{Tobias Holicki, Carsten W. Scherer and Sebastian Trimpe
\thanks{This project has been funded in part by Deutsche Forschungsgemeinschaft (DFG, German Research Foundation) under Germany's Excellence Strategy -EXC 2075 -390740016 and in part by the Cyber Valley Initiative, which is gratefully acknowledged by the authors.}
\thanks{ Tobias Holicki and Carsten W. Scherer are with the Department of Mathematics, University of Stuttgart, Pfaffenwaldring 5a, 70569 Stuttgart, Germany
	(email: \{tobias.holicki, carsten.scherer\}@imng.uni-stuttgart.de)}%
\thanks{Sebastian Trimpe is with the Institute for Data Science in Mechanical Engineering, RWTH Aachen University, Germany, and also with the Intelligent Control Systems Group, Max Planck Institute for Intelligent Systems, Stuttgart, Germany (email: trimpe@dsme.rwth-aachen.de)}%
}

\begin{document}

\maketitle

\thispagestyle{empty}
\pagestyle{empty}

\begin{abstract}
	
For a partially unknown linear systems, we present a systematic control design approach based on generated data from measurements of closed-loop experiments with suitable test controllers.
These experiments are used to improve the achieved performance and to reduce the uncertainty about the unknown parts of the system. This is achieved through a parametrization of auspicious controllers with convex relaxation techniques from robust control, which guarantees that their implementation on the unknown plant is safe.
This approach permits to systematically incorporate available prior knowledge about the system by employing the framework of linear fractional representations.

\end{abstract}


\begin{textblock}{13.1}(1, 15.25)
	\fbox{
		\begin{minipage}{\textwidth}
					\footnotesize \textcopyright 2020 IEEE. Personal use of this material is permitted. Permission from IEEE must be obtained for all other uses, in any current or future media, including reprinting/republishing this material for advertising or promotional purposes, creating new collective works, for resale or redistribution to servers or lists, or reuse of any copyrighted component of this work in other works.
			DOI: \href{https://doi.org/10.1109/LCSYS.2020.3004506}{10.1109/LCSYS.2020.3004506}
		\end{minipage}
}	
\end{textblock}

\section{INTRODUCTION}

\IEEEPARstart{R}{ecently}, learning and data-based control design approaches have received a lot of attention even for linear systems \cite{FerUme20, BocMat18, FazGe19, ZhaHu20,  BerSch16, MarHen17}.
These approaches can often be subsumed under the broad framework of reinforcement learning \cite{KobBag13}, but are still rather diverse \cite{MatPro19}.
In \cite{FerUme20} robust control is combined with a dual design strategy that is used for exploring the closed-loop behavior, while
\cite{BocMat18} employs the system level synthesis framework with an identification step followed by a robust design and an end-to-end analysis.
The approaches in \cite{FazGe19, ZhaHu20} are based on policy gradient methods, while
\cite{BerSch16, MarHen17} rely on Bayesian optimization strategies involving Gaussian processes for tuning the controller parameters. The latter strategies turned out to be very efficient for various applications, in particular, in robotics \cite{DriEng17, CalSey15, RohTri18}.

Bayesian optimization and other direct sampling methods aim to synthesize optimal controllers based on measurements of a closed-loop cost function involving an unknown system $P_0$ to which suitable test controllers are applied \cite{BerSch16, MarHen17, DriEng17, CalSey15, RohTri18}. While these methods have successfully been used in practice, several aspects are subject to current research:
\begin{itemize}[leftmargin=*]
	\item A critical issue is \emph{safety} which means here (and in contrast to the many other interpretations as, e.g., in \cite{DeaTu19}) that the implemented controllers are guaranteed to stabilize the unknown plant $P_0$ \cite{BerSch16, TurKra19}. 
	Such guarantees are not often provided in learning control, which might lead to catastrophic outcomes due to closed-loop instability during the tuning process. To this end, a safe threshold on the cost is introduced in \cite{BerSch16} as an indicator for stability, while \cite{TurKra19} incorporates a robustness objective in terms of classical delay and gain margins.
	\item The choice of a suitable \emph{parametrization} of test controllers is another important issue which aims to keep the evaluations of the cost small even if the set of admissible controllers is large \cite{MarHen16, RobMan11, BanCal17}.
	In \cite{RobMan11}, several naive parametrizations are illustrated and one based on the Youla parametrization is studied. In \cite{MarHen16}, the controller candidates are parametrized in terms of the weights in an LQ design for a given nominal system.
	\item Different ways to incorporate \emph{prior knowledge} is another topic of tremendous importance in these approaches \cite{MarHen16, MarHen17, BocMat18, KobBag13, BerSch15, RohTri18}.
	A linearization of the underlying nonlinear system is used in \cite{MarHen16} for the construction of a parametrization.
	In \cite{MarHen17}, prior knowledge is used for the design of specialized kernels that outperform standard ones, while
	\cite{RohTri18} discusses how to choose hyperparameters from some simulation model.
\end{itemize}
\vspace{-1ex}

In this paper, we propose a systematic parametrization of controllers based on modeling, analysis and design techniques from robust control that can be used for controller tuning/sampling and addresses all of the above concerns at the same time.

We assume that $P_0$ is only partially unknown and employ the linear fractional representation (LFR) framework in order to separate known from unknown (or difficult) components.
Such representations are well-established and flexible modeling tools in robust control \cite{ZhoDoy96, Sch01a}, but they are not often used in learning control.
In particular, LFRs allow for expressing $P_0$ as feedback interconnection of a known linear system $P$ and some unknown or uncertain component $\Del\in\Delf$; the set $\Delf$ captures, e.g., crude guesses on parameter ranges. \emph{Prior knowledge} is thus encoded in the choices of $P$ and $\Delf$.
Dedicated robust design techniques then allow the synthesis of controllers that stabilize the uncertain interconnection and, hence, are guaranteed to stabilize the unknown $P_0$; these techniques ensure \emph{safety}.
In this initial work, we assume that the uncertain component is parametric and construct a \emph{parametrization} based on a partition of the set $\Delf= \bigcup \Delf_k$. The main idea is to use controllers as obtained from a robust multi-objective design problem with guaranteed stability and performance on $\Delf$ and $\Delf_k$, respectively.

\vspace{0.25ex}
\noindent\textbf{Outline.} %
%
The remainder of the paper is organized as follows. After a short paragraph on notation,
we specify the considered learning control problem and discuss its essential ingredients.
Next we propose a systematic parametrization of robust controllers for safely and exploratively evaluating the underlying closed-loop cost function.
We elaborate on the properties of this parametrization and demonstrate its benefits on some numerical examples.

\vspace{0.25ex}
\noindent\textbf{Notation.} %
%
We use the star product ``$\star$'' and all rules for linear fractional transformations (LFTs) as in \cite[Chapter 10]{ZhoDoy96}.
Objects that can be inferred by symmetry or are not relevant are indicated by ``$\bullet$''.

\section{SETTING}

\subsection{Problem Formulation}

We assume that we are given an unknown real system $P_0$ described as
\begin{equation}
\mat{c}{e \\ y}
= P_0
\mat{c}{d \\ u}.
\label{sy0}
\end{equation}
Here $e$ is the controlled output and $d$ is a generalized disturbance (both used to formulate performance specifications), while $y$ is the measurement output and $u$ the control input.
The underlying control problem is to find a controller
\begin{equation}
u = F y
\label{F}
\end{equation}
such that the corresponding closed-loop system, which is referred to as $P_0\star F$, is stable and such that a closed-loop cost function $J$, which encodes the performance specifications, is minimized.
Since $P_0$ is unknown, we aim to find such a controller based on evaluations of the cost function $J$. This amounts to the selection of suitable test controllers, their implementation on the real system $P_0$ and the evaluation of their achieved closed-loop performance.

This is the setting in \cite{BerSch16, MarHen17, DriEng17, CalSey15, RohTri18}, where the individual approaches differ, e.g., in the choice of cost, the available measurements from the plant, the assumed prior knowledge about the plant and the employed controller parametrization.


We confine the discussion to continuous-time linear time-invariant (LTI) systems $P_0$ and the design of state-feedback controller gains $F$ which motivates to choose $y$ as the state $x$ of $P_0$. Moreover, we choose the $H_\infty$-norm cost function
\begin{equation}
J: F \mapsto J(F):=\|P_0 \star F\|_\infty,
\label{DC::eq::org_cost}
\end{equation}
for which ample motivations are found in the robust control literature. Data-based techniques for estimating $H_\infty$-norms have been proposed, e.g., in \cite{RalFor17}. To simplify the exposition, we assume that the measurements of the cost are exact, although it is possible to extend our framework to noisy ones.

\subsection{Encoding Prior Knowledge}

We consider the case that $P_0$ is partially unknown. To this end, we adopt the framework of LFRs \cite{ZhoDoy96} and describe $P_0$ as the interconnection of some known system $P$ given by
\begin{equation}\label{sy}
\mat{c}{\dot x(t) \\\hline z(t) \\ e(t) \\ y(t)}
= \mat{c|ccc}{A & B_1 & B_2 & B_3 \\ \hline
	C_1 & D_{11} & D_{12} & D_{13} \\
	C_2 & D_{21} & D_{22} & D_{23} \\
	I   & 0      & 0      & 0 }
\mat{c}{x(t) \\ \hline w(t) \\ d(t) \\ u(t)}
\end{equation}
in feedback with some unknown part or uncertainty
\begin{equation}\label{De0}
w(t) = \Delta_0z(t).
\end{equation}
Here $\Delta_0$ is a real matrix of suitable dimension and $w$, $z$ are the interconnection variables. Then \r{sy0} admits a state-space representation with the matrix
\begin{equation*}
\arraycolsep=.9pt
\mat{c|cc}{A & B_2 & B_3 \\ \hline
	C_2 & D_{22} & D_{23} \\I&0&0}
\!+\! \mat{c}{B_1 \\ \hline D_{21}\\0}\!\Del_0 (I - D_{11}\Del_0)^{-1} \mat{c|cc}{C_1 & D_{12} & D_{13}}.
\end{equation*}	
By slightly abusing the notation, the latter matrix and the system \r{sy0} are denoted as $\Delta_0\star P$.
Note that the controlled unknown system, the interconnection of \r{sy}, \r{De0} and \r{F}, is then given by
$P_0\star F=(\Delta_0\star P)\star F=\Delta_0\star P\star F$.
Such representations are known to be highly flexible since they permit to effectively capture structural dependencies of models on
uncertain scalar parameters or matrix sub-blocks, which are typically collected on the diagonal of the (structured) uncertainty $\Delta_0$. As an extra advantage, LFRs allow a seamless generalization to multiple heterogeneous (i.e., a mixture of time-varying, non-linear or infinite dimensional) uncertainties collected in a nonlinear feedback operator $w(t)=\Delta_0(t,z)(t)$, but this is not pursued here.

Instead, we adopt the point-of-view that $\Delta_0$ describes unknown parts of the system $P_0$, while the known parts (as, e.g., resulting from first-principle modeling) are captured by $P$. Moreover, we assume that $\Delta_0$ is contained in some known set $\Delf$ of matrices  that is compact and typically given by
\begin{equation*}
\{\diag(\del_1 I, \dots, \del_{m_r} I, \Del_1, \dots, \Del_{m_{\rm f}}) : |\del_j|\leq 1,~ \|\Del_j\|\leq 1\}
\end{equation*}
with (repeated) diagonal and full unstructured blocks on the diagonal, all bounded in norm by one.
As an extreme case, this description does capture the models in \cite{BerSch15, MarHen17, BocMat18} and the ones in \cite{FerUme20, MatPro19, FazGe19, ZhaHu20}, in which it is assumed that nothing aside from linearity is known about $P_0$ and where $\Delta_0$ is just one large unstructured uncertain matrix.

Note that the development of modern robust control has been substantially motivated by the fact that facing completely unknown systems is often not realistic.
By now LFRs are used in tandem with dedicated analysis and design tools from robust control, such as structural singular values or integral quadratic constraints (IQCs) \cite{MegRan97}, which permit to accurately exploit the fine structure of the unknown $\Del_0$.

Therefore, in view of their modeling power, LFRs provide an ideal setting to incorporate prior structural knowledge about a system (through $P$) with unknown to-be-learnt components (through the elements of $\Delta_0$).

\subsection{Safety}

Clearly, guaranteeing stability is a critical issue in learning based approaches since probing the system with gains that are not stabilizing can lead to catastrophes. In contrast to many other approaches as, e.g., in \cite{DriEng17, CalSey15, RohTri18} and aligned with \cite{MarHen17}, we propose to only select controllers that are guaranteed to be robustly stabilizing, i.e., that are taken from the set
\begin{equation}
\F(\Delf) := \{F:~ F \text{ stabilizies }\Del \star P \text{ for all }\Del \in \Delf \}.
\label{DC::eq::set_of_stab_controllers}
\end{equation}
This set is typically much smaller than the set of controllers that are merely required to stabilize $\Delta_0\star P$. However, since $\Delta_0\in\Delf$ is unknown and since we can only rely on the prior knowledge about $\Delf$, there is no other choice than to pick gains from $\F(\Delf)$ in order to ensure a safe operation of the system in closed-loop. The minimal value of $J(F)$ over the set $\F(\Delf)$ is related to the cost of interest as
\begin{equation}
\inf_{F \text{ stabilizes }P_0} J(F) \leq \inf_{F \in \F(\Delf)}J(F),
\label{DC::eq::safety_ineq}
\tag{\ensuremath{\mathcal{S}}}
\end{equation}
in which the gap reflects the price to-be-paid for safety.

\subsection{Motivation for Controller Parametrizations}

For optimizing the cost $J$ it is highly beneficial and an often seen strategy to parameterize a family of test controllers by a few parameters \emph{before} applying an optimization algorithm, especially if the ambient space of controller gains has a large dimension  \cite{MarHen16, RobMan11, BanCal17}.
Formally, such a parametrization is a mapping $\Fc$ with a domain $\mathrm{dom}(\Fc)$ that is contained in a low dimensional ambient space and chosen in order to render the gap in the inequality
\begin{equation}\label{DC::eq::parametrization_gap}
\inf_{F \text{ stabilizes }P_0} J(F)
\leq \inf_{\theta \in \mathrm{dom}(\Fc)}J(\Fc(\theta))
\tag{\ensuremath{\mathcal{P'}}}
\end{equation}
as small as possible. Then, the idea is to minimize the surrogate cost $J \circ \Fc$ over $\mathrm{dom}(\Fc)$ instead of determining the minimum of the original cost $J$.
Since the former minimization problem is formulated in a low dimensional space, it is expected to require substantially fewer evaluations of the cost function for its (approximate) solution.
The gap in \r{DC::eq::parametrization_gap} constitutes the price to-be-paid for this reduction of complexity and is rarely analyzed in the literature.
We stress that it is instrumental to choose a parametrization $\Fc$ such that its values are contained in $\F(\Delf)$ for reasons of safety. Then the gap in \r{DC::eq::parametrization_gap} can even be more precisely identified as the sum of that in \eqref{DC::eq::safety_ineq} and the one in
\begin{equation}\label{DC::eq::parametrization_gap2}
\inf_{F \in \F(\Delf)}J(F)
\leq \inf_{\theta \in \mathrm{dom}(\Fc)}J(\Fc(\theta)).
\tag{\ensuremath{\mathcal{P}}}
\end{equation}

\subsection{Main Contributions}

For some index set $\I$, we propose a novel parametrization $\Fc:\I \to \F(\Delf)$ of auspicious robustly stabilizing controllers based on a partition of $\Delf$. 
It features an a priori safety guarantee without the need to ensure this property through the employed optimization algorithm as proposed, e.g., in \cite{ZhaHu20}.
For its construction, we use advanced robust control techniques that explicitly take the available prior knowledge into account. Based on this parametrization, we show how experimental controller probing allows for controlling the size of the gap in \eqref{DC::eq::parametrization_gap2} by varying the coarseness of the partition of $\Delf$, and how to even reduce the gap in \eqref{DC::eq::safety_ineq} by systematically decreasing the size of $\Delf$ without endangering safety.

In comparison to a standard robust design, which does not utilize data from closed-loop experiments, our approach naturally generates safe controllers with improved performance on the real plant $P_0$.

\section{PARAMETRIZATION OF TEST CONTROLLERS}

\subsection{Construction of the Controller Parametrization}

Let us choose the index set $\I:=\{1, \dots, N\}$ and subsets $\Delf_1, \dots, \Delf_N$ of the uncertainty set $\Delf$ that form the partition
\begin{equation*}
\Delf = {\textstyle \bigcup\limits_{k\in\I}} \Delf_k
\text{ ~with~ }
\mathrm{int} \Delf_k \cap \mathrm{int}\Delf_l = \emptyset
\text{ ~for all~ }
k,l\in\I.
\end{equation*}
With this partition, we construct $\Fc:\I\to \F(\Delf)$ based on the rationale to render $J(\Fc(k))$ for at least one index $k\in\I$ as small as possible, since this leads to the best possible reduction of the gap in \eqref{DC::eq::parametrization_gap2}.
The proposed parametrization assigns to $k\in\I$ a controller $F\in \F(\Delf)$ which reduces
${\sup_{\Del \in \Delf_k}\|\Del \star P \star F\|_\infty}$ 
as much as possible. This means that we are facing a robust multi-objective synthesis problem involving robust stability w.r.t. $\Del \in \Delf$ and worst-case $H_\infty$ performance w.r.t. $\Del \in \Delf_k$. Such problems are usually nonconvex as well as nonsmooth and thus hard to solve systematically.
Still, it is possible to compute good  upper bounds on the corresponding optimal value by solving a linear SDP if relying on so-called multiplier relaxations in robust control. One such relaxation is given in Theorem~\ref{DC::lem::part} and requires to specify a set $\P(\Delf)$ of real symmetric matrices with an LMI description such that
\begin{equation*}
\mat{c}{-\Del^T \\ I}^T P \mat{c}{-\Del^T \\ I} \cl 0
\text{ ~for all~ }
\Del \in \Delf,\ P \in \P(\Delf).
\end{equation*}
We also assume that such multiplier classes $\P(\Delf_k)$ are available for the partition members $\Delf_k$ and for $k=1,\ldots,N$.
A more detailed discussion with concrete choices for such multiplier sets can be found in \cite{Sch05, SchWei00}.

\begin{theorem}
	\label{DC::lem::part} For fixed $k\in \I$ consider the system of LMIs
	\begin{subequations}
		\label{DC::lem::eq::lmi_part}
		\begin{equation}
		\arraycolsep=3pt
		Y\cg 0,\ \ (\bullet)^T \mat{cc|c}{0 & I \\ I & 0 \\ \hline && P}\hspace{-1ex} \mat{cc}{I & 0 \\ -\Ab^T & -\Cb_1^T \\ \hline 0 & I \\ -B_1^T & -D_{11}^T}  \cg 0,
		\label{DC::lem::eq::lmi_partb}
		\end{equation}
		\begin{equation}
		\arraycolsep=3pt
		(\bullet)^T
		\mat{cc|c|cc}{0 & I & &&\\ I & 0 & &&\\ \hline && P_k & &\\ \hline &&& -\ga^2I&0\\&&& 0&I}
		\mat{ccc}{I & 0 & 0\\ -\Ab^T & -\Cb_1^T & -\Cb_2^T \\ \hline
			0 & I & 0 \\
			-B_1^T & - D_{11}^T & - D_{21}^T \\ \hline
			0 & 0 & I \\
			-B_2^T & - D_{12}^T & - D_{22}^T} \cg 0
		\label{DC::lem::eq::lmi_partc}
		\end{equation}
	\end{subequations}
	in the variables $Y=Y^T$, $P\in \P(\Delf)$, $P_k \in \P(\Delf_k)$, $M$, $\ga$ and with the abbreviations
	\begin{equation*}
	\arraycolsep=2pt
	(\Ab, \Cb_1, \Cb_2) \!:=\! (AY\!+\!B_3 M, C_1Y\!+\!D_{13}M, C_2Y\!+\!D_{23}M).
	\end{equation*}
	If these LMIs are feasible, the controller gain $F := MY^{-1}$ satisfies
	$F\in\F(\Delf)$ and $\sup_{\Del \in \Delf_k}\|\Del \star P \star F\|_\infty < \ga$.
\end{theorem}

\vspace{1ex}
The proof of this result is found in \cite{SchWei00}. It shows that
\begin{equation}\label{gaprel}
\inf_{F\in\F(\Delf)} \sup_{\Delta\in\Delf_k}\|\Delta\star P\star F\|_\infty\leq \ga_*(k)
\end{equation}
is satisfied for $\ga_*(k)\!:=\!\inf\{\ga\in\R:\text{ LMIs \eqref{DC::lem::eq::lmi_part} are feasible}\}$.

All this leads us to the construction of the parametrization $\Fc$ as follows: For some fixed small $\eps>0$ and $\ga_*^\eps(k):=(1+\eps)\ga_*(k)$, we assign to $k\in\I$ some gain $\Fc(k)$ with
\begin{equation}\label{Fcprop}
\Fc(k)\in\F(\Delf)\text{ and }\sup_{\Del \in \Delf_k}\|\Del \star P \star \Fc(k)\|_\infty \leq \ga_*^\eps(k).
\end{equation}
We emphasize that both $\ga_*(k)$ and $\Fc(k)$ can be computed by solving a standard semi-definite program.	Still note that, in general, $\ga_*(k)$ is not attained (no optimal controller exists), which motivates the introduction of $\eps$.

In the sequel, we abbreviate the surrogate cost function resulting from the parametrization $\Fc$ as
$$
L(k):=J(\Fc(k))=\|\Delta_0\star P\star \Fc(k)\|_\infty\text{ ~for~ }k\in\I.
$$
Further, let us note at this point that $\Delta_0\in\Delf_{k}$ for some index $k\in\I$ clearly implies
\begin{equation}
L(k)
\leq \sup_{\Del \in \Delf_{k}}\|\Del \star P \star \Fc(k)\|_\infty
\leq \ga_*^\eps(k).
\label{Fckey}
\end{equation}

\begin{remark}
\label{DC::rema::singletons}
It is routine to adapt Theorem~\ref{DC::lem::part} to a singleton $\Delf_k=\{\Delta\}$ with any $\Delta\in\Delf$. This adaptation no longer requires to choose a multiplier class for $\Delf_k$ which promotes a smaller relaxation gap in \eqref{gaprel}.
It also permits the choice $\I=\Delf$ as a highly useful extreme case in our construction and leads to a parametrization $\Fc$ mapping $\Delf$ into $\F(\Delf)$.
\end{remark}

\begin{remark}
As a key difference between the cost $J$ and its surrogate $J \circ \Fc$, the domain of the former consists of the only implicitly defined set of (robustly) stabilizing controllers, while the latter can be evaluated directly. In particular for $\I=\Delf$ as in Remark \ref{DC::rema::singletons}, $J\circ \Fc$ is simply defined on $\Delf$.
\end{remark}

\subsection{Application of the Controller Parametrization}\label{DC::sec::appli}

After having introduced the controller parametrization, the conceptual algorithm of this paper reads as follows. For each $k\in\I$, we can implement the controller $\Fc(k)$ on the system, since it is assured to be stabilizing for $P_0$, and measure the  cost $L(k)$.
A mere minimization over $k\in\I$ then leads to an optimal controller $\Fc(k_\ast)$, and inequality \eqref{DC::eq::parametrization_gap} now reads as
\begin{equation}\label{opt}
\inf_{F \text{ stabilizes }P_0} J(F)\leq L(k_*)\text{ ~for~ }k_*\in\argmin_{k\in\I}L(k).
\end{equation}
Fine partitions of $\Delf$ lead to large index sets $\I$. Instead of considering all $k\in \I$, we can take fewer (random) samples $\{k_1, \ldots, k_M\}$ of $\I$ and obtain a (rough) approximation $\min_{j=1,\dots, M} L(k_j)$ of $L(k_\ast)$.
In particular for the partition with $\I=\Delf$ as described in Remark \ref{DC::rema::singletons}, one can directly employ a whole variety of smarter (derivative free) sampling and optimization strategies, such as Bayesian optimization involving Gaussian processes discussed in \cite{BerSch16, MarHen17}.

\vspace{1ex}

Instead of considering a single (fine) partition, one can as well start from a coarse partition of $\Delf$ and propose adaptive refinement strategies which generate a sequence of controller parametrizations as follows.
Given $\Delf = \bigcup_{k = 1}^N \Delf_k$, determine an index $k^{0} \in \argmin_{j=1,\dots,N} L(j)$. Then generate
a partition of $\Delf_{k^0}$ denoted as $\bigcup_{j = 1}^N \Delf_{k^0_j}$ in order to obtain a refined partition of the original set as
\begin{equation}
\Delf = \Bigl({\textstyle \bigcup\limits_{j = 1,\dots, N}} \Delf_{k^0_j}\Bigr) ~ \cup ~ \Bigl({\textstyle \bigcup\limits_{k = 1,\dots, N,~ k \neq k^0}} \Delf_k\Bigr).
\label{DC::eq::refinement}
\end{equation}
This refined partition yields a new parametrization $\Fc^1$ with corresponding new cost $L^1$ and some next optimal index $k^1 \in \argmin_{j=1,\dots,N} L^1(k^0_j)$. By construction it is guaranteed that $L^1(k^1) \leq L(k^0)$ holds.
This step can be iterated in order to further decrease the value of the surrogate cost.
In Section \ref{DC::sec::algo} we propose a specific algorithm based on this approach which involves, in particular, a concrete strategy for refining given partitions.

\subsection{Reducing the Gap in Inequality \eqref{DC::eq::parametrization_gap2}}

Our setup allows for identifying the sources of the gap in \eqref{DC::eq::parametrization_gap2} and permits to generate systematic refinements towards its reduction.
To illustrate this issue, let us suppose that the relaxation gap in \eqref{gaprel} is small. Then we infer (by the definition of $\Fc$ and for small $\eps>0$) that
$$
\ga_*^\eps(k)\!
\approx \!\inf_{F\in\F(\Delf)}\sup_{\Delta\in\Delf_{k}}\|\Delta\star P\star F\|_\infty
\!\leq\! \sup_{\Delta\in\Delf_{k}}\|\Delta\star P\star \Fc(k)\|_\infty.
$$
On the other hand, for $k\in\I$ with $\Delta_0\in\Delf_{k}$ and if this member $\Delf_{k}$ of the partition is sufficiently small, we have
$$
\sup_{\Delta\in\Delf_{k}}\|\Delta\star P\star \Fc(k)\|_\infty
\approx\|\Delta_0\star P\star \Fc(k)\|_\infty
=L(k).
$$
Hence $\inf_{k\in \I} L(k)$ is close to $\inf_{F\in\F(\Delf)}\|\Delta_0\star P\star F\|_\infty$ which shows that the gap in \eqref{DC::eq::parametrization_gap2} is small.
In conclusion, it is essential that the size of the partition member containing $\Del_0$ and the relaxation gap in \eqref{gaprel} are both small.
Without going into details, we emphasize that the latter can be controlled with the choices of the multiplier sets $\P(\Delf)$ and $\P(\Delf_k)$, through the use of more advanced multi-object control techniques and by applying further refinements in robust control \cite{ArzPea00}, such as incorporating S-variables \cite{EbiPea15} or dynamic instead of static IQCs \cite{MegRan97}.

%

\subsection{Reducing the Gap in Inequality \eqref{DC::eq::safety_ineq}}

Our approach offers the opportunity to even reduce the gap in \eqref{DC::eq::safety_ineq} by identifying a smaller index set $\t\I\subset\I$ with	
\begin{subequations}
	\label{DC::eq::set}
	\begin{equation}
	\Del_0\in \t \Delf := {\textstyle \bigcup_{k \in \t \I}} \Delf_k.
	\end{equation}
	Indeed, this is guaranteed with 	
	\begin{equation}
	\t \I = \{k \in \I : L(k)\leq\ga_*^\eps(k)\}
	\end{equation}
\end{subequations}
since $k \in \I / \t \I$ implies $\Del_0 \notin \Delf_k$ by \eqref{Fckey}.	
Note that $\t \Delf$ can be considerably smaller than the original $\Delf$, which implies that the related set $\F(\t \Delf)$ of robustly stabilizing controllers is (much) larger than $\F(\Delf)$. Thus, replacing $\Delf$ with $\t \Delf$
reduces the cost of safety as expressed by the gap in \eqref{DC::eq::safety_ineq}.

This suggests to repeat our design procedure for $\t \Delf$, which amounts to constructing a new parametrization $\t \Fc$ giving controllers $\t \Fc(k)$ (via Theorem \ref{DC::lem::part}) with which we can perform new closed-loop experiments to evaluate $\t L = J \circ \t \Fc$. The controllers $\t \Fc(k)$ are expected to achieve (considerably) improved closed-loop performance with a smaller gap in \eqref{DC::eq::parametrization_gap}, just due to the reduction of the gap in \eqref{DC::eq::safety_ineq}.
The algorithm proposed in the next section is based on this strategy.

\begin{remark}
	\label{DC::rema::shrinking}
The set $\t \Delf$ is not guaranteed to be convex. Similarly as in \cite{Sch01}, in this case one can express it as union of few convex sets and modify Theorem~\ref{DC::lem::part} by using a robust stabilization objective for each of the individual convex sets;  this purposive design comes along with an increased numerical burden.
\end{remark}

\vspace{0.5ex}

Note that the set $\I$ is constructed based on \eqref{Fckey} which provides an upper bound on the cost $L(k)$ at the index $k$ with $\Del_0 \in \Delf_{k}$. We can also devise a lower bound which can be exploited similarly in order to further reduce $\t \I$ and shrink the gap in \eqref{DC::eq::safety_ineq}.
To this end, observe that standard $H_\infty$ design permits to numerically determine
\begin{equation*}
\ga_\mathrm{nom}(\Del)\! :=\! \inf_{F \text{ stabilizes }\Del \star P}\! \|\Del \star P \star F\|_\infty
\text{ for any fixed }\Del\! \in \!\Delf.
\end{equation*}
Then $\Delta_0\in\Delf_k$ for $k\in\I$ indeed yields the lower bound
\begin{equation*}
\inf_{\Del \in \Delf_{k}}\! \ga_\mathrm{nom}(\Del)
\leq\! \ga_\mathrm{nom}(\Del_0)
\leq \|\Del_0 \star P \star \Fc(k)\|_\infty \leq L(k).
\end{equation*}
Note that this lower bound is not cheap to compute as it involves a numeric minimization of $\ga_\mathrm{nom}$ on $\Delf_k$ for each considered $k \in \I$. In contrast, the upper bound $\ga_\ast^\eps(k)$ is essentially obtained for free while constructing the map $\Fc$.

\section{AN ALGORITHM}\label{DC::sec::algo}

In this section we propose a concrete algorithm that works in higher dimensions \emph{and} aims to exploit \eqref{DC::eq::set}.
It involves the uncertainty box $\Delf = \{\Del(\del): \del_\nu \in \Ic_\nu,\ \nu=1,\ldots,M \}$ where $\Ic_1, \dots,\Ic_M$ are given intervals and in which we use the abbreviation $\Del(\del) := \diag(\del_1 I_{q_1}, \dots, \del_M I_{q_M})$. The related
Algorithm \ref{DC::alg} is motivated by coordinate descent, which currently becomes more popular due to its appearance in machine learning applications.

\SetInd{0.25em}{0.4em}
\begin{algorithm}
	\small
	\SetKwInOut{Input}{input}\SetKwInOut{Output}{output}
	\Input{Number of partitions $N$}
	
	Set $\nu = 1$ and $\Ic_k^p = \Ic_k$ for all $k = 1, \dots, M$\\
	\While{(not terminated)}{
		Choose a uniform partition $\Ic_\nu = \bigcup_{k = 1}^N \t \Ic_k$\\
		Set \hspace{-0.3ex}$\Delf_k \!:=\! \{\Del(\del)\!:\! \del_\nu \!\in\! \t \Ic_k,\, \del_j \!\in\! \Ic_j,\, j \!\neq\! \nu \}$  \hspace{-0.3ex}to  \hspace{-0.3ex}get \hspace{-1ex} $\Delf\! =\! \bigcup_{k = 1}^N\! \Delf_k$\\
		Determine $\Fc(k)$, $\ga_*(k)$ and $L(k)$ for all $k\in \I$\\
		Determine $\t \I$ as in \eqref{DC::eq::set}, set $\Ic_\nu\! =\! \mathrm{convex\, hull}\big(\bigcup_{k \in \t \I} \t \Ic_k\big)$ and update $\Delf$ accordingly\\
		Set $\Ic^p_\nu := \t \Ic_j$ where $j\! \in\! \argmin_{k\in \t \I}L(k)$\\
		Set $\nu = \nu + 1$ if $\nu < M$ and $\nu = 1$ otherwise
	}
	Set $\nu = 1$\\
	\While{(not terminated)}{
		Choose a uniform partition $\Ic_\nu^p = \bigcup_{k = 1}^N \t \Ic_k$\\
		Set $\Delf_k := \{\Del(\del)~:~ \del_\nu \in \t \Ic_k,~ \del_j \in \Ic^p_j, ~ j \neq \nu \}$ to get $\Delf =\big(\bigcup_{k = 1}^N \Delf_k\big) \cup (\bullet)$ as in \eqref{DC::eq::refinement}\\
		Determine $\Fc(k)$, $\ga_*(k)$ and $L(k)$ for all $k\in \{1,\dots,N\}$\\
		Set $\Ic^p_\nu \!:=\! \t \Ic_j$ and $F_\ast \!=\! \Fc(j)$ for $j\hspace{-0.5ex} \in\hspace{-0.5ex} \argmin_{k=1,\dots, N}L(k)$\\
		Set $\nu = \nu + 1$ if $\nu < M$ and $\nu = 1$ otherwise
	}
	\Output{Controller gain $F_\ast$}
	\caption{\small Design via Coordinate-Like Descent}
	\label{DC::alg}
\end{algorithm}

The first loop of the algorithm generates a partition of $\Delf$ by taking a uniform partition only of the interval $\Ic_\nu$ related to the parameter $\delta_\nu$. In line 6, it exploits \eqref{DC::eq::set} in order to shrink
$\bigcup_{k = 1}^N \t \Ic_k$ to a new interval and to generate a reduced parameter set that is guaranteed to contain $\Delta_0$;
moreover, it still has the structure of a hyperrectangle. In line 7 and as input to the second loop, we store those intervals $\Ic_\nu^p$ for which the best performance level is observed. Running this loop $n_1$ times requires to perform $Nn_1$ experiments.

In the second loop, the algorithm adaptively refines those subsets 
of $\Delf$ for which the best closed-loop performance was achieved. This proceeds as in Section \ref{DC::sec::appli}, by
generating sub-partitions along each parameter axis. Again, $n_2$ runs of the loop require $Nn_2$ experimental cost evaluations. In particular, if we let $n_1, n_2 = O(M)$, the number of evaluations grows linearly in the number of unknown parameters $M$ and turns the algorithm applicable even if $M$ is large.


\section{NUMERICAL EXAMPLES}\label{DC::sec::exa}

For numerical illustrations, we consider several modified examples from the library COMPl\textsubscript{e}ib \cite{Lei04} which, unfortunately, does not comprise robust control examples. We let
$\smat{A & B_2 & B_3 \\ C_2 & D_{22} & D_{23}}$ be the matrices $\smat{A & B_1 & B \\ C_1 & D_{11} & D_{12}}$
	in (1.1) from \cite{Lei04} and choose the remaining matrices in order to define $P$ in \eqref{sy} as $D_{11}=0$, $D_{12} = 0$, $D_{21} = 0$,
	\begin{equation*}
	\arraycolsep=2pt
	B_1\!=\!\mat{ccc}{1 & 0 & 1 \\ 0 & 1 & 0 \\ \hdashline \multicolumn{3}{c}{0_{\bullet \times 3}}}\!, ~
	C_1\!=\!\mat{cc:c}{1 & 0  & \\ 0 & 1 & 0_{3\times \bullet} \\ 0 & 1}\text{ and }
	D_{13}\!=\!\mat{c:c}{0  & \\ 0 & 0_{3 \times \bullet} \\ 1 &}\!.
	\end{equation*}
	Further, we take 
	$\Delf \!:=\! \{\diag(\del_1, \del_2, \del_3) : \del_1, \del_2, \del_3  \in [-1, 1] \}$
and suppose that the real system $P_0 = \Del_0 \star P$ is obtained for $\Del_0 = \diag(0.7, -0.1, 0.7)$.
Note that all subsequent algorithms only access the cost
$J(F)=\|P_0\star F\|_\infty$ in \eqref{DC::eq::org_cost} for
chosen gains $F$.
In Theorem \ref{DC::lem::part}, both sets of multipliers $\P(\Delf)$ and $\P(\Delf_k)$ consist of so-called $D/G$-scalings \cite{SchWei00} and we take $\eps:= 0.05$ in \r{Fcprop}. With Theorem \ref{DC::lem::part} for
the trivial partition $N=1$, we can as well compute an upper bound $\ga_{\rm rp}$ for the robust performance synthesis problem as in
\begin{equation}\label{rp}
\inf_{F\in\F(\Delf)}\sup_{\Del\in\Delf}\|\Delta\star P\star F\|\leq\ga_{\rm rp}.
\end{equation}
Any learning based design results in a controller with performance level in between $\ga_{\rm rp}$ and the best achievable nominal performance
$\ga_{\rm nom}:=
\inf_{F \text{ stabilizes }P_0}\! \|P_0 \star F\|_\infty$.

Let us now employ Algorithm \ref{DC::alg} with $N = 6$ and by using $n_1$, $n_2$ iterations in the first and second loop, respectively. The achieved performance levels $\ga^{n_1,n_2}$ for $(n_1,n_2)=(0,6)$ and $(n_1,n_2)=(3,3)$ are depicted in Table~\ref{DC::tab::results}.

We compare these results with a based-line learning approach which aims to minimize the cost $J$ without employing a controller parametrization, similarly as done in \cite{MarHen17, FazGe19, ZhaHu20}.
To this end, we use a deterministic direct search method \cite{AudDen02} which is initialized with a robust controller
as obtained by computing $\ga_{\rm rp}$ in \r{rp}. We rely on the Matlab implementation in \texttt{patternsearch} and
denote by $\ga_\mathrm{ps}^k$ the achieved performance level after $k$ evaluations of the cost function.
Let us emphasize at this point that, for the considered examples, all iterates of \texttt{patternsearch} are stabilizing $P_0$ without any particular precautions; this is in stark contrast to 
other direct optimization algorithms such as \texttt{bayesopt}.
Moreover, we also employ \texttt{patternsearch} for minimizing $J \circ \Fc$ if making use of our parametrization $\Fc$ for
$\I=\Delta$ as described in Remark \ref{DC::rema::singletons}. It is then initialized in $\Del = 0$ and $\ga_{\mathrm{ps}\Fc}^k$ denotes the achieved performance level after $k$ cost evaluations.
%


The results in column $3$ and $4$ of Table \ref{DC::tab::results} demonstrate the benefit of exploiting our controller parametrization over a direct minimization of the cost for only a few (here $k=36$) iterations, despite the gap in \eqref{DC::eq::parametrization_gap}
and even coming along with safety guarantees.
This can be attributed to the fact that the dimension of the ambient controller gain space is larger than $20$ for some examples, since $k=200$ iterations lead to further improvements of performance as shown in column $2$, but without guarantees for stability along the iterations.

The results in columns $5$ and $6$ for Algorithm \ref{DC::alg} show performance levels mostly similar to $\ga_{\mathrm{ps}\Fc}^{36}$ and for an identical number $36$ of evaluations of the cost.
The comparison of $\ga^{3,3}$ with $\ga^{0,6}$ for the second group of examples reveals the benefit of exploiting \eqref{DC::eq::set} in the first loop of Algorithm \ref{DC::alg}.

Let us finally point to the first and last column in Table \ref{DC::tab::results} in order to illustrate
the general benefit of our safe learning approach over a standard robust design, by finding
controllers with (sometimes even drastically) improved closed-loop $H_\infty$ performance for $P_0$, which even comes close
to the optimal level $\ga_\mathrm{nom}$ in some cases.


\begin{table}
	\vspace{1.5ex}
	\caption{\rm Nominal performance and performance achieved by controllers resulting from minimizing $J$ and $L$ via \texttt{patternsearch}, from Algorithm \ref{DC::alg} and from a standard robust design for several modified examples from \cite{Lei04}}
	\label{DC::tab::results}
	\centering
	\setlength{\tabcolsep}{5pt}
	\renewcommand{\arraystretch}{1.2}
	\begin{tabular}{@{}l@{\hskip 3.5ex}rrr@{\hskip 4.5ex}rrr@{\hskip 4.5ex}r@{}}
		\toprule
		Name &$\ga_{\rm nom}$ 
& $\ga_\mathrm{ps}^{200}$ & $\ga_\mathrm{ps}^{36}$ & $\ga_{\mathrm{ps}\Fc}^{36}$ & $\ga^{0, 6}$ &$\ga^{3,3}$ & $\ga_{\rm rp}$
\\ \hline
		AC3 & 3.07 & 4.82 & 5.18 & 3.38 & 3.45 & 3.28 & 5.92 \\
		AC6 & 5.31 & 5.90 & 5.90 & 5.40 & 5.40 & 5.55 & 5.91  \\
		AC11 & 2.72 & 2.80 & 5.55 & 2.73 & 2.79 & 2.75 & 6.57 \\
		HE2 & 1.67 & 1.93 & 7.25 & 4.78 & 4.82 & 4.87 & 7.29 \\
		REA2 & 0.51 & 0.60 & 0.62 & 0.52 & 0.51 & 0.51 & 0.62 \\
		DIS2 & 0.64 & 0.98 & 1.05 & 0.65 & 0.66 & 0.67 & 1.55 \\
		TG1 & 3.52 & 3.94 & 4.04 & 3.60 & 3.64 & 3.64 & 4.35 \\
		ROC6 & 2.40 & 3.23 & 3.53 & 3.36 & 3.36 & 3.36 & 3.84 \\ \hdashline
		AC2 & 0.11 & 0.18 & 0.20 & 0.16 & 0.16 & 0.13 & 0.21\\
		HE4 & 5.30 & 9.06 & 9.17 & 6.77 & 6.50 & 5.32 & 9.31 \\
		DIS1 & 2.42 & 2.65 & 7.40 & 4.93 & 4.97 & 2.67 & 7.40 \\
		MFP & 3.47 & 5.43 & 5.31 & 5.99 & 5.39 & 4.29 & 8.03 \\
		NN4 & 1.00 & 1.14 & 1.79 & 1.52 & 1.53 & 1.02 & 2.57 \\
		NN16 & 0.95 & 1.48 & 1.48 & 1.02 & 1.04 & 0.97 & 1.50 \\
		\bottomrule
	\end{tabular}
\vspace{-4ex}
\end{table}

\section{CONCLUSIONS}

We consider LTI systems affected by an unknown parameter $\Del_0$ contained in some known set $\Delf$ and propose strategies to design controllers based on generated data from measurements of closed-loop experiments with suitable test controllers.
For the systematic selection of auspicious candidates, we propose a new controller parametrization induced by a partition of $\Delf$ and based on advanced robust control techniques.
In particular, this parametrization ensures that all candidates are robustly stabilizing which guarantees that their implementation on the real system is safe.
Interestingly, it even offers the possibility to systematically generate subsets of $\Delf$ which are guaranteed to contain the unknown $\Del_0$.


The concept admits immediate extensions to output feedback control if relying on existing design techniques with robust stability and performance guarantees. It is as well easily possible to consider $H_2$-norm cost criteria on an infinite or finite time-horizon or in discrete-time.
While the employed modeling and design tools from robust control offer much more flexibility in terms of capturing time-varying, dynamic or nonlinear unknown components, the systematic construction of controller parametrization along the presented lines remains largely open in such situations.


\addtolength{\textheight}{-12cm}   



%

%



\bibliographystyle{IEEEtran}
\bibliography{literaturs}

\end{document}